\documentclass[11pt]{article}
\usepackage{amsmath,latexsym}
\usepackage{amsfonts}
\newtheorem{theorem}{Theorem}[section]

\def\slfrac#1#2{\hbox{\kern.1em %
 \raise.5ex\hbox{\the\scriptfont0 #1}\kern-.11em %
 /\kern-.15em\lower.25ex\hbox{\the\scriptfont0 #2}}}

\newcommand{\eqn}[1]{(\ref{#1})}

\newcommand{\eeq}{\end{equation}}
\newcommand{\beql}[1]{\begin{equation}\label{#1}}
\newcommand{\bsq}{{\vrule height .9ex width .8ex depth -.1ex }}

\newcommand{\CC}{{\mathbb C}}

\newcommand{\RR}{{\mathbb R}}

\newcommand{\bJ}{{\bf J}}

\newcommand{\bM}{{\bf M}}

\newcommand{\sD}{{\cal D}}
\newcommand{\sE}{{\cal E}}
\newcommand{\sH}{{\cal H}}

\newcommand{\sM}{{\cal M}}

\makeatletter
\def\@sect#1#2#3#4#5#6[#7]#8{\ifnum #2>\c@secnumdepth
     \def\@svsec{}\else
     \refstepcounter{#1}\edef\@svsec{\csname the#1\endcsname.\hskip .75em }\fi
     \@tempskipa #5\relax
      \ifdim \@tempskipa>\z@
        \begingroup #6\relax
          \@hangfrom{\hskip #3\relax\@svsec}{\interlinepenalty \@M #8\par}%
        \endgroup
       \csname #1mark\endcsname{#7}\addcontentsline
         {toc}{#1}{\ifnum #2>\c@secnumdepth \else
                      \protect\numberline{\csname the#1\endcsname}\fi
                    #7}\else
        \def\@svsechd{#6\hskip #3\@svsec #8\csname #1mark\endcsname
                      {#7}\addcontentsline
                           {toc}{#1}{\ifnum #2>\c@secnumdepth \else
                             \protect\numberline{\csname the#1\endcsname}\fi
                       #7}}\fi
     \@xsect{#5}}
\def\@begintheorem#1#2{\it \trivlist \item[\hskip \labelsep{\bf #1\ #2.}]}

\def\plain{plain}\ifx\fmtname\plain\csname fi\endcsname
     
     \input docstrip
     \preamble

     Do not distribute the stripped version of this file.
     The checksum in the header refers to the documented version.

     \endpreamble
     \generateFile{here.sty}{t}{\from{here.doc}{}}
     \endinput
\fi
\ifcat a\noexpand @\let\next\relax\else\def\next{%
    \documentstyle[here,doc]{article}\MakePercentIgnore}\fi\next
\ifx\@Hxfloat\@Hundef\else\expandafter\endinput\fi
\let\@Hxfloat\@xfloat
\def\@xfloat#1[{\@ifnextchar{H}{\@HHfloat{#1}[}{\@Hxfloat{#1}[}}
\def\@HHfloat#1[H]{%
\expandafter\let\csname end#1\endcsname\end@Hfloat
\vskip\intextsep\vbox\bgroup\def\@captype{#1}\parindent\z@
\ignorespaces}
\def\end@Hfloat{\egroup\vskip \intextsep}

\makeatother

\catcode`\@=11
\renewcommand{\section}{
        \setcounter{equation}{0}
        \@startsection {section}{1}{\z@}{-3.5ex plus -1ex minus
        -.2ex}{2.3ex plus .2ex}{\large\bf}%
        }
\catcode`@=12
\begin{document}
\begin{center}
{\Large {\bf The Schr\"{o}dinger Operator with Morse Potential on the Right Half Line}}\\
\vspace{1.5\baselineskip}
{\em Jeffrey C. Lagarias 
\footnote{This work was  supported by the NSF
under Grants DMS-0500555 and DMS-0801029.  }}\\
\vspace*{.2\baselineskip}
University of Michigan \\
Ann Arbor, MI 48109-1043 \\
{\tt lagarias@umich.edu} \\
\vspace*{2\baselineskip}
(August 10, 2009-revision) \\
\vspace{1.5\baselineskip}
{\bf Abstract}
\end{center}
\noindent
This paper studies the Schr\"{o}dinger operator with Morse
potential $V_{k}(u) = \frac{1}{4} e^{2u} + ke^{u}$
on a right  half-line $[u_0, \infty)$, and determines
the Weyl asymptotics of  eigenvalues for constant boundary
conditions at the endpoint $u_0$. In consequence it  
 obtains information on the location of zeros of the Whittaker
 function $W_{\kappa, \mu}(x)$, for fixed real parameters $\kappa, x$
 with $x>0$, 
 viewed as an entire function of the complex variable  $\mu$. In this case
 all zeros  lie on the imaginary axis, with the possible 
 exception, if $\kappa > 0$, of a finite number of real
zeros  which lie  in the interval $-\kappa <  \mu < \kappa$.
We obtain an asymptotic formula for the number of  zeros 
$N(T) = \{\rho~| ~ W_{\kappa, \rho}(x)=0, ~|Im(\rho)| < T\}$ of the form
$N(T) = \frac{2}{\pi} T \log T + \frac{2}{\pi}(2\log 2 -1 - \log x) T +O(1).$
Parallels are observed with  zeros of the Riemann zeta function.\\

\noindent

%
%
%
\setlength{\baselineskip}{1.2\baselineskip}

\section{Introduction}

The one-dimensional potential function $V_{A,B}(u)=Ae^{2u}+ Be^u$  with $A >0$ and $B$ real
was proposed in 1929 by Philip  Morse \cite{Mo29}
as an approximation to 
the (radial) quantum-mechanical potential for diatomic molecules, cf. Flugge \cite[Problem 70]{Fl74}. 
From a 
quantum physics viewpoint this
corresponds to the Schr\"{o}dinger operator 
$$
-\frac{d^2}{du^2} + V_{A, B}(u),
$$
viewed either on the  line $(-\infty, \infty)$ or on the left  half-line
$(-\infty, u_0]$, with a suitable boundary condition imposed at  the  endpoint $u_0$.
 By rescaling the
independent variable  we may without loss of generality
reduce to the case  $A= \frac{1}{4}$,  and in that case
relabel the other parameter,  setting  $B=k$, obtaining 
the potential $V_k(u) = \frac{1}{4}e^{2u} + k e^u$.
The Morse potential on the line, with $k=0$,
arises elsewhere in physics  in connection with
the Liouville model on the hyperbolic plane (Grosche \cite{Gr88}) and in
describing motion on the hyperbolic plane with a
magnetic field  (Ikeda and Matsumoto \cite{IM99}). The latter authors also relate it to
the Selberg trace formula on compact hyperbolic surfaces.

In both the line and left half-line cases  the Schr\"{o}dinger operator 
has an absolutely continuous spectrum of multiplicity one 
supported on the half-line $[0, \infty)$,  plus a
finite (possibly empty) discrete spectrum, lying strictly
below the absolutely continuous spectrum.
On the  line, the operator is essentially
self-adjoint for the standard domain of smooth rapidly decreasing functions, and is exactly
solvable.
Formulas for its resolvent (Green's function)
have been obtained, e.g. \cite[Prop. 4.1(i)]{IM99}.  
The absolutely continuous spectrum is described
using generalized eigenfunctions given
in terms of  Whittaker $M_{\kappa, \mu}(x)$-functions, and the  discrete spectrum 
has eigenfunctions given using rescaled Laguerre polynomials,
cf. Ismail and Koelink \cite{IK08}. 
%
\subsection{Morse potential on right half-line}

We treat the Morse potential $V_k(u)$ on the right half-line
$[u_0, \infty)$, with constant boundary
conditions imposed at the endpoint $u_0$. The spectrum of this operator is 
 quite different from the cases above:  it is pure discrete, simple, and bounded below.
This follows because the Morse potential is bounded below
and is unbounded as $u \to \infty$. 
The discrete  eigenfunctions on the half-line
are expressible in terms of Whittaker
functions $W_{\kappa, \mu}(x)$, in which $\mu$ is the eigenvalue parameter.
We consider
\beql{103}
\left(-\frac{d^2}{du^2} + V_{k}(u) \right) \psi(u) = E \psi(u),
\eeq
and  in \S2 we observe for all complex $E$
there is a unique eigenfunction (up to scaling) in $L^2([u_0, \infty), du)$
with eigenvalue $E=-\mu^2$ given by 
\beql{103b}
\psi(u, E) := e^{-\frac{1}{2} u} W_{-k,\pm \mu}(e^{u}),
\eeq
The function $\psi(u, E)$  is well-defined since  $W_{k, \mu}(x)= W_{k, -\mu}(x)$  
for all $\mu \in \CC$,   when $k$  is real and 
$x=e^u$ is positive real.
The allowed real eigenvalues $E_0  < E_1 < E_2 < ...$ in the discrete spectrum
are  selected by the boundary conditions at the
left endpoint $u_0$, with corresponding value 
$x_0= e^{u_0}$,
since we are in the limit point case at the singular endpoint $u=+\infty$.

This paper 
determines  the Weyl asymptotics of the spectrum for
fixed constant boundary conditions, obtaining in \S4 that
$$
\# \{ E_n \le T\} = c_1 \sqrt{T} \log T + c_2 \sqrt{T} + O(1)
$$
as $T \to \infty$. 
Dirichlet boundary conditions correspond to zeros of the Whittaker-function
in the $\mu$-variable, holding the other variables fixed. The results above,
applied with these boundary conditions, 
lead to determination of the location and asymptotic distribution of the
zeros of these Whittaker functions, given in \S5 below. Namely, these zeros
all lie on the imaginary axis, with finitely many exceptions,
 with asymptotics given in \eqn{111a} below. The exceptional zeros 
 exist only if  $\kappa=-k >0$, and they all lie on
 the real axis in the open interval $(-|k|, |k|)$.


%
\subsection{Analogy with 
distribution of 
Riemann zeta zeros}

Our main interest in this family of operators arose 
from number theory,   
 to view them as   ``toy models'' for  
operators whose eigenvalues represent
zeros of the Riemann zeta function.
We note there has been a persistent effort to find
a spectral interpretation of the Riemann zeta zeros. In this
regard, by a ``Hilbert-Polya" operator we will mean an (unbounded)
self-adjoint operator on a Hilbert space whose spectrum encodes the
zeta zeros, 
such that the self-adjointness of the operator
encodes the Riemann hypothesis. Here we  consider one-dimensional 
Schr\"{o}dinger operators on a  half-line from this viewpoint.

To compare the distribution of
Morse half-line eigenvalues with zeros of the Riemann zeta function $\zeta(s)$, 
 in our context the $s$-variable in $\zeta(s)$ corresponds to  the $\mu$-variable in the
Whittaker function $W_{\kappa, \mu}(x)$ via 
$s :=\mu - \frac{1}{2}$;
 this  variable change maps the imaginary axis in 
the $\mu$-variable to the critical line $Re(s)= \frac{1}{2}$ in the $s$-variable.
One point of this paper is that the function
\[
Z_{1}(s)  := Z_1(s; \kappa, x_0) = W_{\kappa, s-\frac{1}{2}}(x_0),
\]
for fixed real $\kappa$ and fixed real $x_0>0$  has  a number of  properties
analogous to the  Riemann $\xi$-function
$$\xi(s)= \frac{1}{2} s(s-1)\pi^{-\frac{s}{2}} \Gamma(\frac{s}{2})
\zeta (s).$$
 Namely, in Theorem~\ref{th51} we establish 
the following properties of $Z_1(s)$.\\

(1)  $Z_1(s)$ is an entire function of order $1$ and
maximal type, which is real on the real axis and  is real
on the line   $Re(s) = \frac{1}{2}$. \\

(2) $Z_1(s)$ satisfies the functional equation
$$Z_{1}(s)=Z_{1}(1-s).$$

(3) The number $N_{T}^{\pm}(Z_1(s))$   of zeros $\rho$ of $Z_1(s)$ with $|\rho|<T$ is given
as $T \to \infty$ by
\beql{111a}
N_{T}^{\pm}(Z_1(s)) = \frac{2}{\pi} T \log T + 
\frac{2}{\pi}( 2\log 2 - 1 -\log x_0) T +O(1).
\eeq

(4)  All but finitely many of the zeros
of $Z_1(s)$ lie on the line $Re(s)= \frac{1}{2}$. The finite set of 
``exceptional zeros''  all  lie on the real axis, in the open interval 
$( \frac{1}{2}- |\kappa|,\frac{1}{2}+ |\kappa|)$. There are no exceptional
zeros if $\kappa \le 0$. \\

(5)  All zeros of $Z_1(s)$ are simple zeros, except for a possible
double zero at $s= \frac{1}{2}$.\\

\noindent 
The zeros $\rho= \frac{1}{2}+ i \gamma$ of $Z_1(s)= Z_1(s; \kappa, x_0)$ have the following
spectral interpretation:
  $\rho$ is a zero of $Z_1(s; \kappa, x_0)$
if and only if  $E=  \gamma^2$ is an eigenvalue of the 
Schr\"{o}dinger operator for the Morse potential $V_k(u)$
on the half-line
$[u_0, \infty)$ for $k= - \kappa$ and $x_0=e^{u_0}$, 
with  Dirichlet boundary conditions taken at the left endpoint $u_0$. 
The self-adjointness of this problem leads to  the ``Modified Riemann
Hypothesis" property (4). This spectrum is
bounded below by $E \ge - \kappa^2$, and the eigenvalues $E < 0$
correspond to exceptional zeros of $Z_1(s; \kappa, x_0)$ lying on the real axis.

In comparison with $Z_1(s)$, the  Riemann $\xi$-function $\xi(s)$ has properties (1) and (2). 
Analogous to property (3), the number of
zeros $N_T^{\pm}(\xi(s))$ of the Riemann $\xi$-function with $|\rho|<T$  satisfies
\beql{112}
N_T^{\pm}(\xi(s)) = \frac{1}{\pi} T\log T +
\frac{1}{\pi} (-\log 2\pi -1) T + O(\log T),
\eeq
cf. Edwards \cite[Sec. 6.7]{Ed74} or Titchmarsh \cite[Theorem 9.4]{TH86}.
Analogues of properties (4) and (5) are conjectured to hold for $\xi(s)$.
Property (4) is 
a ``modified Riemann Hypothesis",  which permits 
a finite number of exceptional
zeros to occur on the real axis. The Riemann
$\xi$-function has no exceptional zeros,
but the possibility of exceptional zeros remains open  for various generalizations of the
Riemann zeta function,  e. g. automorphic $L$-functions, cf. Iwaniec and Sarnak \cite{IS00}.
Concerning property (5), the zeta function 
 is reported to have simple zeros  for the first $10^{13}$ zeros on the critical line.

A possible weakness of the Morse potential  as  a``toy model"   for a spectral interpretation
of  the zeta zeros is that  its
eigenvalue distribution only  reproduces the main terms in
the asymptotic number of zeros, it  does not correctly model   the  ``lower order"
asymptotic behavior of zeta zeros. The detailed structure of scaled zeta zeros is
conjectured to be described by the GUE distribution of random matrix theory,
and the functions here do not reproduce this behavior. This detailed behavior
presumably encodes all the number theory having to do with the distribution
of prime numbers. This is discussed further in \S5. A putative ``Hilbert-Polya"
operator Schr\"{o}dinger potential, if it exists, will likely have exotic features to produce
such features in its spectrum.

%
\subsection{Schrodinger operators and de Branges spaces as spectral models for zeta zeros}

Recently we observed  that, provided the Riemann hypothesis
holds,  there is a natural candidate for  a ``Hilbert-Polya" 
operator  within  the framework of the de Branges theory
of Hilbert spaces of 
entire functions  (\cite{La06}, \cite{La07}). 
The de Branges  theory provides  model operators  for a certain class
of operators, which  includes  generalizations of operators 
from singular Sturm-Liouville problems of the type considered here.
This theory  gives a representation of
the associated operator as a (self-adjoint) {\em canonical system} of differential equations
(defined in \S5),  in which
each zeta zero $\rho= \frac{1}{2}+ i \gamma$ corresponds to a 
simple real eigenvalue $\gamma$ of
the resulting canonical system, so that the operator has a simple spectrum 
that is  unbounded above and below. 
The results in
((\cite{La06}, \cite{La07}) show that if appropriate Riemann hypotheses hold then
 the associated canonical systems {\em must exist.}  We call these
 Hilbert-Polya canonical systems.

The connection to Schr\"{o}dinger operators is that, in  
favorable circumstances,  a canonical system may be  formally  transformed
to a Schr\"{o}dinger operator on a half-line. This transform 
squares the eigenvalues, so that the resulting Schr\"{o}dinger spectrum is bounded below.
 If this transformation could legitimately be done for
the putative ``Hilbert-Polya" de Branges spaces above,
then  it would result in a ``Hilbert-Polya"
operator that is a Schr\"{o}dinger operator on a half-line.
If so, it would 
be very interesting to determine the  associated potential of this  operator. 
Note that this  ``Hilbert-Polya" potential, if it exists at all, may be highly singular,  involve
fractal features, and be describable only as a distribution. We remark that
additional symmetries of the problem indicate that this
Schr\"{o}dinger  operator   necessarily would impose   Dirichlet
boundary conditions  at the left endpoint.

The  Schr\"{o}dinger operators treated in this paper 
may provide  hints about features
of this  possible ``Hilbert-Polya"  potential. 
In particular, they  give the dependence 
of eigenvalue asymptotics on  the endpoint of the interval. 
Furthermore, even  if such a ``Hilbert-Polya" Schr\"{o}dinger operator does not exist, these
examples  may still give hints about the canonical system of the de Branges space
associated to RH. Namely, 
there exists a de Branges space with explicitly known canonical
system for which the nonlinear transform above can be carried out to
give a  ``toy model" Morse potential Schr\"{o}dinger operator  with 
$k = -\frac{1}{2}$ (\cite{La08}). The form of this particular canonical system then gives
hints about features of the associated canonical system coefficients
of the Hilbert-Polya canonical system,
which is guaranteed to  exist if  RH holds.

%
\subsection{Previous work}

There  has been  a long history of work on special
 functions whose zeros  mimic 
 the zeros of the Riemann zeta function, some of them associated
 to differential equations.  Polya \cite{Po26}
 observed in  1926  that the zeros of the $K$-Bessel function $K_{\mu}(x)$ for  fixed
$x>0$ lie on the imaginary
axis and have asymptotics similar to that of the Riemann zeta function.   
In fact Polya  actually studies a function given by an integral representation, which 
in retrospect was shown to be  a $K$-Bessel function, cf. Gasper \cite{Ga08}. 
In 1927 Polya \cite{Po27} gave
 a second proof of the property that these zeros lie on the imaginary
 axis which made use of  a second-order differential
equation in an auxiliary variable. Another 
``toy model"  discussed by Polya is  described 
in Titchmarsh \cite[Sect. 10.1]{TH86}. It  involves a shifted sum of two $K$-Bessel functions,
discussed more recently in   Gasper \cite{Ga08}.
The present work can be viewed as obtaining an extension to a one-parameter family of
Polya's results, in the sense that a  $K$-Bessel function is essentially a special
case of the Whittaker functions considered here, namely
$K_{\mu}(w) = \sqrt{\frac{\pi}{2w}} W_{0, \mu}(2w).$ 
Recent work of Biane \cite{Bi08} 
notes an operator-theoretic
parallel between the $K$-Bessel function and the Riemann $\xi$-function,
and makes a very interesting new connection between these two functions in terms 
of probability theory.

Aside from the ``toy model" aspect,
our results provide some new information on the location of zeros of Whittaker functions. 
Here there has been considerable previous work locating  zeros of $W_{\kappa, \mu}(z)$,
in connection with various applied problems.
This work 
mostly considers the case of fixed
$\kappa$ and $\mu$, and varies the (complex) variable
$z$.  In 1915 Milne \cite{Mi15} studied real $\kappa$ and $\mu$ 
and determined the nature of zeros in the $z$-variable on the positive real axis
in special cases. 
 In 1941 Tvestkoff \cite{Tv41} announced results on the location
of all zeros with  $\kappa$  and $\mu$ real and  $z=x >0$ a positive
real, and in 1950 Tricomi \cite{Tr50} gave
detailed proofs.  Their work determines  the locations of the 
finite set of  ``exceptional zeros" in (4) above. 
In 1960 Dikii \cite{Di60} and  Dyson \cite{Dy60} independently  considered the 
 case of $\kappa$ real,
$\mu$ purely imaginary, and  $z$ a (nonzero) complex variable.
Here $W_{\kappa, \mu}(z)$ is multi-valued 
in the $z$-variable with a singularity at $z=0$, and  the $z$-plane is cut 
along the nonpositive real axis. Their work arose from 
analysis of the stability of a hydrodynamic model of an incompressible atmosphere
having density decreasing exponentially with height, with wind velocity increasing linearly
with height,  see Dikii \cite{Di60b} and Case \cite{Ca60}, repectively.
Both Dikii and Dyson showed  that  all the zeros in the complex variable $z$ 
in the region $-\pi \le arg(z) \le \pi$   lie on the positive real axis, 
and that there are infinitely many such zeros, having the singular point $z=0$ as a limit point. 
 The results here complement these analyses by determining the location and density of  zeros when
$\kappa$ is real and $z$ is positive real.

Our  interest in de Branges spaces
was influenced by recent 
work of Burnol (\cite{Bu02}, \cite{Bu03}, \cite{Bu06}), who
described a scattering associated to the Fourier and
Hankel transforms, and  pointed out a de Branges space interpretation.   
Burnol constructed  certain Dirac operators and Schr\"{o}dinger operators 
(corresponding to  particular canonical systems) which have
eigenvalue asymptotics  of the general form $cT \log T + O(T)$.
In these papers Burnol did not explicitly note these asymptotics, 
but he elsewhere derived similar
eigenvalue asymptotics  for related operators  (\cite{Bu04}).

%
\subsection{Methods}

The proofs  given here are entirely classical, within the framework of Titchmarsh \cite{Ti62}.
They  break no new ground methodologically.
The  main point  of this paper is  rather to explain 
an analogy of these operators with a ``Hilbert-Polya"
operator for the Riemann zeta  zeros. 
In particular, we  explicitly determine the 
eigenvalue asymptotics \eqn{111a}, with all details
included, for the convenience of number theorists. \\

%
\subsection{Contents of Paper}

The  detailed contents of this paper are as follows.

In  \S2 we determine the 
$L^2$-eigenfunctions of  the Schr\"{o}dinger operator
with the Morse potential on the half-line $[u_0, \infty)$, and
determine its solutions for constant  boundary conditions
at the left endpoint $u_0$,
particularly Dirichlet boundary
conditions (Theorem \ref{th21}). This is a singular boundary value problem.
Using Sturm-Liouville theory we deduce certain
monotonicity properties of the eigenvalues as functions
of the parameters (Theorem \ref{th22}). 

In \S3 we show that  the density of
eigenvalues $\#\{ E \le T\}$  of the Schr\"{o}dinger operator is asymptotic to
$c_1 \sqrt{T} \log T + c_2 \sqrt{T} + O(1)$, for certain constants
$c_1, c_2$ (Theorem ~\ref{th34}).

In \S4, we specialize to  Dirichlet boundary
conditions,  and apply these results to Whittaker functions $W_{\kappa, \mu}(x)$,
viewed as a function of the complex variable $\mu$.
We establish that the zeros of $Z(\mu) := W_{\kappa, \mu}(x)$
 lie on the imaginary $\mu$-axis, with finitely many exceptions, which
themselves lie on the real axis, and we determine  their
asymptotic density to height $T$ to be given by \eqn{111a} above (Theorem ~\ref{th51}).

In \S5 we discuss  analogies  of these results with zeta functions in 
number theory and with de Branges spaces.
In this connection, it
is well known that one can relate de Branges spaces to Schr\"{o}dinger operators
(Sturm-Liouville operators) in various ways, see 
Dym \cite{Dy71}, Dym and McKean \cite{DM76}, and
Remling \cite{Re02}.   Specific
de Branges spaces having number theory content were
related   to Schr\"{o}dinger and Dirac operators
by  Burnol \cite{Bu01}, \cite{Bu02}, \cite{Bu03}, \cite{Bu06}, \cite{Bu06b}.

In \S6 we make some concluding remarks, discussing
analytic spectral invariants associated to these operators,
the   Weyl-Titchmarsh $m$-function of a Schr\"{o}dinger operator and
an analogous quantity in the de Branges theory, which we call
the de Branges $m$-function.

There are two appendices. 
 Appendix A reviews  properties of Whittaker 
functions $M_{\kappa, \mu}(z)$ and 
$W_{\kappa, \mu}(z)$ as functions of three complex variables $(\kappa, \mu, z)$.
Appendix B gives formulas for  the principal Weyl-Titchmarsh $m$-functions associated
to  Morse potentials $V_k(u)$ on the right half line. This function encodes spectral
data about the potential, and we show it is  a ratio of Whittaker $W$-functions with shifted
$\kappa$-parameter values, shifted by a constant.  It is well known that a Schr\"{o}dinger
equation potential $V(u)$ can be uniquely reconstructed from the principal Weyl-Titchmarsh
$m$-function under general conditions.

\paragraph{Notation.}
We follow the mathematical convention that Hilbert spaces have scalar  products
$<f,g>$ that are linear in the first factor and conjugate-linear 
in the second factor, as in  Coddington and Levinson \cite{CL55} and
Levitan and Sargsan \cite{LS75}. This is done for  compatibility with the theory of 
de Branges  \cite[p. 50]{deB68} discussed in \S5. References to Reed and Simon \cite{RS75},
who use the physics convention that Hilbert spaces are conjugate-linear
in the first factor, must be adjusted for this fact. 
We set $L^2(u_0, \infty):= L^2[ (u_0, \infty), du]$.
The variables $a, b, c, d, k, t, u, v, x, y$ denote real variables.
The variables $z= x+iy, s= \sigma+it, \kappa= a+ib, \mu=c+id$
denote complex variables.

\paragraph{Acknowledgments.}
 Christian  Remling supplied 
helpful comments and references, and  
J.-F. Burnol made
useful comments  on his work and
on $m$-functions. I thank  them
and the reviewers for helpful comments
and corrections.

%
%
%

\section{Morse Potential on the Right Half Line}

Let $k$ be an arbitrary real number.
We consider the Schr\"{o}dinger operator
$H = - \frac{d^2}{du^2} + V(u)$ for the
Morse potential 
$$
V_{k}(u) = \frac{1}{4} e^{2u} + k e^{u}
$$
on the half-line $[u_0, \infty)$. 
This is equivalent to treating the more general
potential
$$
V_{A, B}(v) = A e^{2v} + Be^{v}
$$
with $A >0$ and $B$ real on the half-line $[v_0, \infty)$, since
by a translation of the time variable $v = u - \frac{1}{2}\log \frac{A}{4}$,
we reduce to a potential of the form $V_{k}(u)$ above
with $k= \frac{2B} {\sqrt{A}}$
on the half-line $[u_0, \infty)$ with $u_0= v_0 - \frac{1}{2}\log \frac{A}{4}$.

%
%
%

\begin{theorem}~\label{th21}
The  Schr\"{o}dinger equation 
\beql{303}
\left(-\frac{d^2}{du^2} + V_{k}(u) \right) \psi(u) = E \psi(u)
\eeq
for the Morse potential  $V_{k}(u) = \frac{1}{4} e^{2u} + k e^{u}$
on a half-line $[u_0, \infty)$ has the following properties.

(1) The two-dimensional complex vector space $\sE_{E}( V_{k}(u))$
of solutions $\varphi(u)$ is given by
\beql{304a} 
\varphi(u) = e^{-\frac{u}{2}}f(e^u)
\eeq
in which $f(x)$ is any solution to Whittaker's differential equation
\beql{311}
\left(\frac{d^2}{dx^2} +  (  - \frac{1}{4} + \frac{\kappa}{x}
+ \frac{\frac{1}{4} - {\mu}^2}{x^2})\right)f(x) =0,
\eeq
with parameters $(\kappa, \mu)$ given by $\kappa =-k$ and $\mu^2 = -E.$

(2) 
For arbitrary $E \in \CC$  there is a one-dimensional
subspace of $\sE_{E}(V_{k}(u))$ of solutions that
belong to  $L^2(u_0, \infty)$. Setting $E=z^2$,
this  subspace is spanned
by a solution  $\psi(u, E)$ given by
\beql{305}
\psi(u, z^2) := e^{-\frac{u}{2}} W_{-k, iz}(e^{u}) = 
e^{-\frac{u}{2}} W_{-k, -iz}(e^{u}),  
\eeq
where $W_{\kappa,\mu}(x)$ denotes the  Whittaker $W$-function.

(3) The potential $V_{k}(u)$  is of limit point type at $+\infty$.
If the constant boundary conditions 
\beql{306}
(\cos \alpha) \psi(u_0) + (\sin \alpha) \psi^{'} (u_0) = 0.
\eeq
for fixed $0 \le \alpha < 2\pi$
are imposed at the left endpoint  $u=u_0$, 
then the Schr\"{o}dinger
equation is self-adjoint. This operator has a pure discrete simple
real spectrum, bounded below. For Dirichlet 
boundary conditions for parameter values $k \ge 0$, the  smallest eigenvalue
$E_0 >  0$; for parameter values  $k < 0$, one has
$E_0 > -k^2.$ 
\end{theorem}

\paragraph{Remark.}
The  case $\alpha=0$ (resp. $\alpha=\frac{\pi}{2}$)
gives Dirichlet (resp. Neumann)
boundary conditions at the left endpoint.
Recall that for general boundary conditions  \eqn{306}, there is
 no uniform lower bound on the spectrum as $\alpha$
is varied over
$0 \le \alpha \le 2 \pi$. 
The eigenvalue $E_0=E_0(\alpha)$ is a continuous
function of $\alpha$ except for  a jump discontinuity at $\alpha=0$.

\paragraph{Proof.} 
(1) Whittaker's differential equation \eqn{311} with parameters
 $(\kappa, \mu)$ was studied in 1904 by Whittaker \cite{Wh04}.
Properties of its solutions are given in 
Whittaker and Watson \cite[Chap. 16]{WW63}, 
Buchholz \cite{Bu69}. Note that the parameters $\pm \mu$ give
the same solution space, since $\mu$ only enters the
differential equation via $\mu^2$.  For
fixed $(\kappa, \mu)$ \eqn{311}  has 
a two-dimensional complex vector space of solutions $\sE_{E}( V_{k}(u))$
(for $z$ on the positive real axis).
This vector space is spanned by the 
Whittaker functions $(W_{\kappa, \mu}(z), M_{\kappa, \mu}(z))$,
except when $\kappa$ is a nonpositive integer.
(Another basis is  $(W_{\kappa, -\mu}(z), M_{\kappa, -\mu}(z))$. Note
that $W_{\kappa, \mu}(z)= W_{\kappa, -\mu}(z)$ but 
 in general $M_{\kappa, \mu}(z) \ne M_{\kappa, -\mu}(z)$.
One calculates by direct substitution that if $f(z)$ is any solution
to Whittaker's equation \eqn{311},
then $\varphi(u) := e^{-\frac{u}{2}} f(e^{u})$ satisfies the
Schr\"{o}dinger equation
\beql{312}
\left( -\frac{d^2}{du^2} + V_{-\kappa}(u) \right) \phi(u)= -\mu^2 \phi(u).
\eeq
Choosing  $(\kappa, \mu)=(-k, \pm \sqrt{-E})$,
these functions span
the  two-dimensional space of solutions $\sE_{E}( V_{-\kappa}(u))$
to the Schr\"{o}dinger equation \eqn{303}.

(2) The  Whittaker function 
$W_{\kappa,\mu}(z)$ (\cite[Chap. 16]{WW63}) is 
a solution $f(z)$ of Whittaker's equation \eqn{311}
which has  the property of rapid decay as $z= x \to \infty$ 
(on the real axis), made unique by the normalization of its 
asymptotics
\beql{313}
 W_{\kappa, \mu}(x) = 
e^{-\frac{1}{2} x} x^{\kappa}( 1 + O\left(\frac{1}{x}\right) )
\eeq
as $x \to \infty$. 
Viewed in the complex domain,
the same asymptotics holds for $|arg(z)| \le \pi - \epsilon < \pi$,
with the implied constant in the error term depending 
on $\epsilon > 0$. 
All other  linearly independent solutions to Whittaker's differential equation
 increase exponentially as $x \to \infty$.
Choosing $\kappa= -k$, the
 particular choice $\varphi(u):=\psi(u,z^2)$ 
coming from the Whittaker function
is the only solution (up to a scalar
multiple) that satisfies  $\varphi(u) \in L^2(u_0, \infty)$ for (any)
finite $u_0$; it has double-exponential decay as $u \to \infty$,
roughly like $e^{-\frac{1}{2} e^u}$, as specified by the
asymptotic formula \eqn{313}. All other solutions have
double-exponential growth as $u \to \infty$. The solutions \eqn{305} thus
 parametrize the subordinate solutions of the Schr\"{o}dinger operator
with Morse potential on the right half-line, in the sense of Gilbert and Pearson \cite{GP87}.

(3) The property that $V_{k} (u)$ is positive for sufficiently large $u$
 puts the Schr\"{o}dinger equation 
in the limit point
case at $u= \infty$ (cf. Reed and Simon \cite[Theorem X.8]{RS75}).
In particular for each real $E$ there is at most  one solution
to the Schr\"{o}dinger equation with eigenvalue $E$
(up to a scalar multiple) that belongs to 
$L^2([u_0, \infty); du)$. Here such a solution was exhibited
for all  complex $E$, given above.

It is known
that for  continuous potentials $V(u)$ on $[u_0, \infty)$ 
that are bounded below, 
the Schr\"{o}dinger operator $-\frac{d^2}{du^2}+V(u)$ with
(separated) boundary conditions given by 
\beql{314a}
(\cos \alpha) f(u_0)  + (\sin  \alpha) f^{'}(u_0)=0
\eeq
for fixed $0 \le \alpha < 2\pi$ 
is  self-adjoint in $L^2([u_0, \infty); du).$ 
To describe the explicit
domain of this operator, let 
$AC[u_0, \infty]$ denote the set of absolutely continuous functions
$f(u)$ on $[u_0, \infty)$ whose (almost everywhere defined) derivative
is in $L^2[u_0, \infty)$. The associated dense domain for
self-adjointness is 
\begin{eqnarray*}
\sD_{\alpha}  & := & \{ f(u) \in C^1[u_0, \infty)\cap L^2([u_0, \infty);du)~|~
 f'(u) \in AC[u_0, \infty], \\
&& \left(-\frac{d^2}{du^2} + V(u) \right)f(u) \in
L^2([u_0, \infty);du), 
~\mbox{and} ~\eqn{314a}~ \mbox{holds} \},
\end{eqnarray*}
cf. Reed and Simon \cite[p. 144]{RS75}.

The fact that $V_{k}(u)$ is continuous, with 
$V_k(u) \to \infty$ as $u \to \infty$ is sufficient to imply
that it has pure discrete spectrum, bounded below,
for the boundary conditions
\eqn{306}, for any choice of $\alpha$ (Levitan and Sargsan \cite[Sect. 4.1,
Lemma 1.2]{LS75}.) This spectrum is  simple
because for each real $E$ there is at most one
solution in $L^2([u_0, \infty); du)$.
The spectrum is bounded below,
for any fixed constant boundary conditions, because the potential 
$V_{k}(u)$ is
bounded below.

The Sturm-Liouville theory applies to real potentials having
a pure discrete simple spectrum bounded below. It implies that
the  eigenvalues
can be numbered $E_0 < E_1< E_2 < ...$ 
and the corresponding 
eigenfunctions $\psi_0< \psi_1<\psi_2 <...$ can be chosen real-valued.
The Sturm-Liouville
theory shows  that the eigenfunction $\psi_n(u)$
has exactly $n$ sign changes
 on the domain $[u_0, \infty)$, so has $n$ zeros; the bottom
eigenfunction $\psi_0(u)$ has no changes of sign and can
be chosen positive. 
(See Atkinson \cite[Chap. 8]{At64}), 
and Coddington and Levinson \cite[Chap. 8]{CL55}.)

We now impose Dirichlet boundary conditions $\alpha=0$,
and lower bound the minimum eigenvalue $E_0$
using the Sturm comparison
theorem. For $k \ge 0$ the eigenvalues of $V_k(u)$
are bounded below by those of $V_{0}(u)$ 
(with same boundary conditions on the endpoint $u_0$), 
since $V_{k}(u) \ge V_0(u)$ on the entire line. 
It suffices to show that the minimum Dirichlet eigenvalue $E \ge -k^2$
for $V_{k}(u)$ with $k < 0$, and $E_0>0$ when $k =0$.
The function $V_{k}(u)$ for $k< 0$ has global 
minimum value $-k^2$,  which occurs at the unique point  $u= \log (-2k),$
while for $k=0$,   $V_{k}(u)$ is strictly positive.
 We use the Sturm comparison theorem for
 the potential
$\tilde{V}_{A, u_0, E}(u) := A (u- u_0)^2-E$ for 
$A >0$, so that 
$-\frac{d^2}{du^2} + V_{A, u_0}(u)$ on $[u_0, \infty)$
is just a translated
and  scaled version of the harmonic oscillator potential 
on $[0, \infty)$, with eigenvalues shifted
by $E$.  For a given $\epsilon >0$ ,  we take $E= k^2 +\epsilon$
and then can choose $A>0$ just small enough that
$$
\tilde{V}_{A, u_0, E}(u) = A(u- u_0)^2-k^2-\epsilon \le
V_{k}(u), ~~~u_0 \le u < \infty.
$$

The Sturm comparison theorem in the singular case
(see Coddington and Levinson \cite[Chap. 8, Theorem 1]{CL55}) 
now shows that the 
minimum eigenvalue $E_0$ of
$V_{k}(u)$ has
\beql{381a}
E_0 \ge \tilde{E}_0 - k^2 - \epsilon,
\eeq
where $\tilde{E}_0$ is the minimum eigenvalue of 
$\tilde{V}_{A, 0,0}(u)= A u^2$
on $[0, \infty)$ with the same boundary conditions at $u=u_0$.
It is well known that the oscillator Hamiltonian $-\frac{d^2}{du^2} +u^2$
on the half-line $[0, \infty)$ with Dirichlet conditions at $u=0$ has a
strictly positive spectrum with minimum eigenvalue  $\tilde{\tilde{E}}_0=3$.
Under rescaling of $A$ we retain the property $\tilde{E_0} > 0$,
and the inequality $E_0 >  -k^2$ then follows from \eqn{381a}
on letting $\epsilon$ become sufficiently small. Finally, for  $k=0$ we may omit $\epsilon,$
and so obtain $E_0>0$. 
$~~\bsq$ \\

 The Sturm theory implies 
monotonicity properties of the eigenvalues as
the parameters $k, u_0$ are varied, with boundary conditions
remaining fixed.

%
%
%

\begin{theorem}~\label{th22}
Let $V_k(u)= \frac{1}{4} e^{2u} + k e^{u}$.
Let $E_{n, \alpha}(u_{0}, k)$ denote the $n$-th eigenvalue 
$(n \ge 0)$ from
the bottom of the spectrum of the operator
$H= -\frac{d^{2}}{du^{2}}+ V_{k}(u)$ on the half-line $[u_{0}, \infty)$
with boundary conditions
\beql{381}
(\cos \alpha)  \psi(u_0) + (\sin \alpha)\psi^{'}(u_0)=0, 
\eeq
with $\alpha$ being fixed.  For each $n \ge 0$ the following hold. \\

(i) For fixed parameter $k \ge 0 $, the $n$-th eigenvalue $E_{n, \alpha}(u_{0}, k)$ 
is a strictly increasing function of $u_0$.\\

(ii) For fixed parameter $k <0$, the $n$-th eigenvalue $E_{n, \alpha}(u_{0}, k)$ 
is  a strictly increasing function of $u_0$, when $u_{0} \ge \log 2|k|$.\\

(iii) For variable $k$, with fixed $u_{0}$,
the $n$-th eigenvalue $E_{n, \alpha}(u_{0}, k)$ is a strictly increasing 
function of $k$.
\end{theorem}

\paragraph{Proof.}
By a change of parameter $v=u-u_0$  we shift the left hand endpoint to $v=0$.
Then each of the assertions (i)-(iii) above can be viewed as varying the potential,
while keeping the left endpoint and the boundary condition fixed. In all cases the
potentials are bounded below, and the variation increases the potential pointwise.
The results then follow from comparison of
the minimax characterization of the $n$-th eigenvalue.
(Reed and Simon \cite[XIII.1]{RS78}).
$~~~~\Box$

%
%
%

\section{Asymptotic Formula for Eigenvalue Density} 

We obtain estimates for the number of
eigenvalues using a version of  Weyl asymptotics, which
says that  the number of eigenvalues $E \le T$ of a Hamiltonian
$H= p^2 + V(q)$ with $V(q)$ bounded below can be 
approximated by the volume in the phase space of
the region $H \le T$. 

The half-line condition requires that we restrict
the $q$-variable in the phase space  to $[u_{0}, \infty)$
in estimating the phase space volume. 
The following result is a slight modification of a result of 
Titchmarsh \cite[Theorem 7.4]{Ti62})
giving a  rigorous  estimate of this kind.

%
%
%

\begin{theorem}~\label{th33}
Let $V(u)$ be a real-valued $C^1$-potential on $[u_0, \infty)$
which for all sufficiently large $u$ is increasing
and convex downwards. Then the Schr\"{o}dinger operator
 $- \frac{d^2}{du^2} + V(u)$  on $[u_0, \infty)$
with boundary condition \\
$$(\cos \alpha) \psi(u_0) + (\sin \alpha) \psi'(u_0)=0,$$
has pure discrete simple spectrum. 
Let $N(T; \alpha, u_0)$ count the number of eigenvalues $E_n < T$
for this boundary condition. 
Then there is a constant $T_1$ such that  for all  $T\ge T_1$ there is a unique solution
$u_T \ge u_0$ to $V(u, T)= T$, and  
\beql{351}
N(T; \alpha, u_0) = \frac{1}{\pi}\int_{u_0}^{u_T}\sqrt{T - V(u)} du + O(1)
\eeq
holds for $T_1 \le T< \infty$.
\end{theorem}

\paragraph{Proof.} Suppose that $V(u)$ is increasing 
and convex downwards for $u \ge u_1$.
The convexity
condition implies for $u \ge u_1$ that  $V(u)$ is strictly increasing 
and increases at least at a linear rate.
Thus  $V(u) \to \infty$ as $u \to \infty$, so 
the Schr\"{o}dinger operator with given boundary conditions
has a pure discrete simple spectrum.

Titchmarsh \cite[Theorem 7.4]{Ti62} proves the
estimate  \eqn{351} assuming as hypotheses that $V(u)$
is continuous,  increasing 
and convex downward on 
the whole interval $[u_0, \infty)$, i.e. that one can take
$u_1 = u_0$. A slight modification of his
argument handles the case above.
We set $T_1 := \max\{ V(u): u_0 \le u \le u_1\}$
For $T >  T_1$, the strict increasing property of $V(u)$
for $u \ge u_1$ implies that  
the equation $V(u)=T$ has a unique solution, so 
that the 
set of values $\{ u: V(u) \le T\}$ is an interval $[u_0, u_T]$.
Titchmarsh's argument estimates
 the number of oscillations
in the real-valued $L^2$ -solution 
$\psi(u; E)$ with eigenvalue $E$ on $[u_0, \infty)$
(not necessarily satisfying the boundary  conditions).
We replace his integral estimate $0 \le I_2(Y) \le \frac{2}{3}$
with $I_2(Y)= O(1)$; the extra $C^1$-condition
on $V(u)$ on the interval $[u_0, u_1]$ is imposed to guarantee that  
the relevant integrals exist. $~~~\bsq$ \\
 
We apply this result to the Morse potential. Note that the parameter $k$
in the Morse potential has no influence on the 
first two terms in the resulting asymptotic formula  \eqn{361} below.

%
%
%

\begin{theorem}~\label{th34}
For the Morse potential $V_k(u) = \frac{1}{4} e^{2u} + k e^u$
on $[u_0, \infty)$ with boundary conditions
$(\cos \alpha) \psi(u_0) + (\sin \alpha) \psi'(u_0)=0,$
the eigenvalue density satisfies 
\beql{361}
N(T; \alpha, u_0) = \frac{1}{\pi} \sqrt{T} \log \sqrt{T} + 
\frac{1}{\pi}\left(2\log 2 - 1 - u_0\right) \sqrt{T} + O(1).
\eeq
as $T \to \infty$, and the implied constant in the O-symbol
depends on  both $k$ and $u_{0}$.
\end{theorem}

\paragraph{Proof.} 
The Morse potential $V_{k}(u)$ is increasing and convex downwards
for all sufficiently large $u$, so Theorem~\ref{th33}
applies to give 
\beql{363}
N(T, \alpha) = \frac{1}{\pi}\int_{u_0}^{u_T}\sqrt{T - V_k(u)} du + O(1).
\eeq
as $T \to \infty$. For  $k \ge 0$ the Morse potential is
increasing and convex downward on the entire real line. However when $k<0$
there is a critical value $u_1$ to the left of which this no longer  holds. 

It remains to estimate the integral in \eqn{363}, with
$V_{k}(u) = \frac{1}{4}e^{2u}+ke^{u}$.
We first suppose that $u_{0}$ is such that $V_{k}(u)$ is 
increasing on $[u_{0}, \infty)$, which holds for all $u_{0}$ if $k \ge 0$
and for $u_{0} \ge u_{1} := \log 2|k|$ if $k <0$. We assume in what
follows that $u_{0} \ge \log (6|k|+1)$, and treat the case
when $u_{0} < \log (6|k|+1)$ at the end of the proof.
We suppose $T$ large enough to be in the region where $V_{k}(u)=T$ has
a unique solution 
and  introduce the positive constant $T^{\ast}$ given by
$(T^{\ast})^{2}= T +k^{2}$ .
We evaluate the
integral
\beql{364}
I := \int_{u_0}^{u_T}\sqrt{T - V_k(u)} du
=\int_{u_0}^{u_T} 
\sqrt{ (T^{\ast})^{2}-\left( \frac{1}{2}e^{u}+k\right)^{2} }du.
\eeq
We change variable $u$ to $w$ determined by
$$
w= \frac{T^{\ast}}{\frac{1}{2}e^{u}+k}.
$$
Since $V_{k}(u)$ is monotone increasing on $[u_{0}, \infty)$ 
the change of variable is single-valued,
and on the integration interval  $w$ decreases from 
$w_0= \frac{T^{\ast}}{\frac{1}{2}e^{u_0}+k}$ to $w_{T}=1$. 
We have 
$$
dw=- \frac{T^{\ast}}{(\frac{1}{2}e^u +k)^2}\frac{1}{2}e^{u}du
= -w(1-\frac{kw}{T^{\ast}}) du,
$$
and obtain
\begin{eqnarray*}
I &= &-\int_{w_{0}}^{1} 
\sqrt{ (T^{\ast})^{2}-\frac{(T^{\ast})^{2}} {w^{2}}}
\left(\frac{1}{1-\frac{kw}{T^{\ast}}}\right)
\frac{dw}{w} \\
&= & T^{{\ast}}\int_{1}^{w_{0}}\sqrt{1-\frac{1}{w^{2}}}
\left(\frac{1}{1-\frac{kw}{T^{\ast}}}\right)
\frac{dw}{w} 
\end{eqnarray*}
in which  
$\frac{T^{\ast}}{w_{0}}=
\frac{1}{2}e^{u_{0}}+k.$
We let $w=e^{v}$, and obtain
\beql{366}
I= T^{\ast}\int_{0}^{v_{0}}\sqrt{ 1-e^{-2v} }
\left(\frac{1}{1-\frac{ke^{v}}{T^{\ast}}} \right)dv,
\eeq
in which 
\beql{367}
v_{0}= \log w_{0}= \log T^{\ast}- \log(\frac{1}{2}e^{u_{0}}+k),
\eeq
where by hypothesis $\frac{1}{2}e^{u_{0}}+k > |k|+\frac{1}{2}$.
Now
$$
e^{v}=w=\frac{T^{\ast}}{\frac{1}{2}e^{u}+k}
$$
so
$$
T^{\ast}= \left(\frac{1}{2}e^{u} + k\right)e^{v}.
$$
The hypothesis on $u_{0} > \log (6|k|+1)$when $k<0$ ensures that 
\beql{368}
0 \le \frac{|k|e^{v}}{T^{\ast}} < \frac{1}{2}.
\eeq
Thus we may legally  expand the integrand of \eqn{366} in a power series
in $ \frac{ke^{v}}{T^{\ast}}$ and obtain $I=M+R$ with
\begin{eqnarray*}
M & := & T^{\ast}\left(\int_{0}^{v_{0}} \sqrt{1-e^{-2v}}dv\right) \\
R & := & T^{\ast}\left(
\sum_{j=1}^{\infty}\left( \frac{k}{T^{\ast}}\right)^{j}
\int_{0}^{v_{0}}e^{jv}\sqrt{1-e^{-2v}}dv \right).
\end{eqnarray*}
Now we have
\begin{eqnarray*}
M & = & T^{\ast}v_{0} + 
T^{\ast}\left(\int_{0}^{v_{0}}(\sqrt{1-e^{-2v}} - 1)dv\right)\\
&=& T^{\ast}\left(\log T^{\ast}- \log (\frac{1}{2}e^{u_{0}}+k) \right) \\
&& +T^{\ast}\left( \int_{0}^{\infty}(\sqrt{1-e^{-2v}}-1)dv + 
O(\frac{1}{(T^{\ast})^2}) \right) 
\end{eqnarray*}
in which the constant in the $O$-symbol depends on $k$ and $u_{0}$.
In the last line we used the estimate
$$
|\int_{v_{0}}^{\infty}(\sqrt{1-e^{-2v}}-1)dv |\le
\int_{\log T^{\ast}- c_{0}}^{\infty}\frac{1}{2}e^{-2v}dv \le 
O\left(\frac{1}{(T^{\ast})^{2}}\right).
$$
Now we use the identity
$$
 \int_{0}^{\infty}(1- \sqrt{1-e^{-2v}})dv= 1- \log 2,
$$
to  obtain
\beql{318a}
M = T^{\ast}\log T^{\ast} + \left(\log2 - 1- \log (\frac{1}{2}e^{u_{0}}+k)
\right) T^{\ast} + O(1),
\eeq
and the constant in the $O$-symbol depends on $k$ and $u_{0}$.
The remaining term is split as $R=R_{1}+R_{2}$, in
which
\begin{eqnarray*}
R_{1} &:= & T^{\ast}\left(\sum_{j=1}^{\infty} 
\left(   \frac{k}{ T^{\ast}}\right)^{j} \int_{0}^{v_{0}} e^{jv}dv\right), \\
R_{2} & := & T^{\ast}\left(\sum_{j=1}^{\infty}
\left(  \frac{k}{ T^{\ast}}\right)^{j}
\int_{0}^{v_{0}}\left(\sqrt{1-e^{-2v}} -1\right) e^{jv}dv\right).
\end{eqnarray*}
The first term above is
\begin{eqnarray*}
R_{1} &=& T^{\ast}\left(\sum_{j=1}^{\infty}(\frac{k}{T^{\ast}})^{j}
(\frac{1}{j}(e^{jv_{0}}-1))\right) \\
&=& T^{\ast}\left( -\log (1 - \frac{ke^{v_{0}}}{T^{\ast}})
 +\log(1-\frac{k}{T^{\ast}})\right)\\
 &=& T^{\ast}\left(-\log( 1 - \frac{k}{\frac{1}{2}e^{u_{0}} +k})+
O(\frac{1}{T^{\ast}}) \right) \\
  &=& T^{\ast}\left( 
  \log(\frac{1}{2}e^{u_{0}}+k) - \log(\frac{1}{2}e^{u_{0}})
  \right) +O\left(1 \right).
\end{eqnarray*}
The second term is estimated by
\begin{eqnarray*}
|R_{2}| & \le & T^{\ast} 
\sum_{j=1}^{\infty}\left(\frac{|k|}{T^{\ast}} \right)^{j}
\int_{0}^{v_{0}}e^{(j-2)v}dv \\
&\le & |k|(1-e^{-v_0}) + \frac{|k|^2}{T^{\ast}} v_0+ 
T^{\ast} \sum_{j=3}^{\infty}
\left(\frac{|k|}{T^{\ast}} \right)^{j}
\left( \frac{T^{\ast}}{\frac{1}{2}e^{u_{0}}+k}\right)^{j-2} \\
&\le & O\left(1\right) + \frac{k^2}{T^{\ast}}\sum_{j=3}^{\infty}
\left(\frac{|k|}{\frac{1}{2}e^{u_{0}}+k}
\right)^{j-2} =O\left(1\right ).
\end{eqnarray*}
Putting these estimates together gives
$$
R = R_1 + R_2 =T^{{\ast}}\left( 
\log (\frac{1}{2}e^{u_{0}}+k) -\log (\frac{1}{2}e^{u_{0}}) \right) + O(1),
$$
and
\beql{369}
I = T^{{\ast}}\log T^{\ast} +\left(2\log2 - 1 -u_{0}\right)T^{\ast}
+O\left(1\right),
\eeq
provided that  $u_{0}\ge \log (6|k|+1)$.

In the remaining case where $u_{0} \le \log (6|k|+1)$,
 we split the integral into a portion from $u_{{0}}$ to 
$u_{1}= \log (6|k|+1)$,
and from $[u_{{1}}, \infty)$.  The integral from $[u_{1}, \infty)$
is estimated by \eqn{369} while that on $[u_{0},u_{1}]$ by
\beql{371}
\int_{u_{0}}^{u_{1}}\sqrt{T- V_{k}(u)}du =
\int_{u_{0}}^{u_{1}}\sqrt{T -O(1)}du = (u_{1}-u_{0})\sqrt{T}
+O\left(\frac{1}{\sqrt{T}}\right),
\eeq
in which the implied constant in the O-symbol depends on $u_{0}$ and $k$.
Note that as $T ^{\ast} \to \infty$, we have
\beql{372}
T^{\ast} = \sqrt{T+k^{2}}= \sqrt{T}\left(1+ \frac{k^{2}}{2T}+
 O(\frac{1}{T^{2}}) \right)= \sqrt{T} +
O\left(\frac{1}{\sqrt{T}}\right).
\eeq
Combining this with \eqn{371} shows that \eqn{369} remains valid for
all intervals $[u_{0},\infty)$. Next, substituting\eqn{372} in \eqn{369}
yields
\beql{373}
I= \sqrt{T} \log \sqrt{T} +
\left(2\log 2 - 1-u_{0} \right)\sqrt{T} + O(1).
\eeq
Substituting this estimate for \eqn{364}  in \eqn{363}  gives the desired estimate \eqn{361} of
$N(T, \alpha)$.
$~~~\bsq$

%
%
%

\section{Zeros of Whittaker Functions in  the $\mu$-Variable}

For Dirichlet boundary conditions  the  eigenvalues
of the Morse potential on the
half-line  are
specified by Theorem~\ref{th21}
in terms of the zeros of Whittaker functions.
We combine this with the results of \S4 to deduce
the distribution of zeros of these Whittaker functions
in the $\mu$-variable. The function $Z_1(s)$ in the introduction
is given as $Z_1(s)= Z(s+ \frac{1}{2})$, where 
$Z(\mu):= W_{\kappa, \mu}(e^{u_0})$ below.
 
%
%
%

\begin{theorem}~\label{th51}
For fixed real parameters $\kappa$ and $u_0$
the   Whittaker function 
\beql{511}
Z(\mu) := W_{\kappa, \mu}(e^{u_0})
\eeq
is an entire function of $\mu$ of order $1$ and maximal type.
It is real-valued on both the real axis and imaginary axis. 
It is an even function of $\mu$, i.e. it satisfies
the ``functional equation''
\beql{512}
Z(\mu)= Z(-\mu).
\eeq

(1) The zeros of $Z(\mu)$  lie on the real and imaginary
axes.  All zeros of $Z(\mu)$ are simple zeros, except
for a possible zero at $\mu=0$, which if it occurs will be a
double zero.

(2) For each $\kappa$ the function $Z(\mu)$ has finitely many real zeros.
There are no real zeros when $\kappa \le 0$, and for $\kappa \ge 0$ they
all lie in the open
interval $-\kappa < \mu < \kappa$.

(3) There are infinitely many  imaginary zeros.
The total number $N(Z(\mu))$
 of zeros with $|Im(\mu)| \le T$
satisfies the asymptotic formula
\beql{513}
N(Z(\mu)) = \frac{2}{\pi} T \log T + 
\frac{2}{\pi}( 2\log 2 - 1 -u_0) T +O(1),
\eeq
as $T \to \infty$, where the $O$-constant depends on $\kappa$.

\end{theorem}

\paragraph{Proof.} We set $x=e^{u}$. It is well known
that for fixed $\kappa$ and $x$ the Whittaker function $W_{\kappa,\mu}(x)$
is an entire function of $\mu$, see Appendix A.  
The function $W_{\kappa, - \mu}(x)$ satisfies Whittaker's
differential equation with the same parameters $(\kappa, \mu)$, and
it also has the same asymptotics as $x \to \infty$, namely
$$
 W_{\kappa, \mu}(x) = e^{-\frac{1}{2} x} x^{\kappa}
( 1 + O\left(\frac{1}{x}\right) ).
$$
It follows that 
\beql{514}
W_{\kappa, -\mu}(x) = W_{\kappa, \mu}(x)
\eeq
holds for all parameters $(\kappa, \mu)$ and positive real $x$.
(It then holds universally under analytic continuation in
the $x$-variable, see  Appendix A.) The functional
equation \eqn{512} follows, taking $x=e^{u_0}$.  This also implies 
that the function $G(\mu) := Z(\sqrt{\mu})$ is an entire function of $\mu$.

We next establish the real symmetry when
$\kappa$ is real-valued
and $x=e^{u_0}$ is positive-real valued (on the principal branch
of the Whittaker function). For $2\mu$ not an integer
we have the representation (\cite[Sec. 16.41]{WW63})
\beql{516} 
W_{\kappa, \mu}(x) =
\frac{\Gamma(-2\mu)}{\Gamma(\frac{1}{2} -\kappa-\mu)}
M_{\kappa, \mu}(x) + 
\frac{\Gamma(2\mu)}{\Gamma(\frac{1}{2} -\kappa+\mu)}
M_{\kappa, -\mu}(x).
\eeq
Now we assert that, for real $\kappa$ and positive real $x$,
applying complex conjugation to \eqn{516} yields
\beql{516a}
W_{\kappa, \bar{\mu}}(x) = \overline{W_{\kappa, \mu}(x)}.
\eeq 
To justify \eqn{516a} , we note that  for positive real $x$, on the principal branch of
the Whittaker $M$-function we have 
$$
M_{\bar{\kappa}, \bar{\mu}}(x) = \overline{M_{\kappa, \mu}(x)}
$$
(see \eqn{905} in the Appendix) and that 
$\Gamma(\bar{s})= \overline{\Gamma(s)}$ for $s\in \CC,$ since
it is real on the real axis. This justifies \eqn{516a}, and 
in fact for real $\kappa$ one has $\overline{W_{\kappa, \mu}(z)}=W_{\kappa, \bar{\mu}}(\bar{z})$
for complex $z$ with $|\arg(z)|< \pi$.
We conclude using \eqn{516} that for real $\kappa$ and positive real $x$,
the function $W_{\kappa, \mu}(x)$ is real when $\mu$ is on the real axis.
In addition, under the same assumptions, the conjugation symmetry  \eqn{516} 
together with the symmetry \eqn{514} implies that
$W_{\kappa, \mu}(x)$ is real when $\mu=it$ is on the imaginary axis.

For real $\kappa$ and positive real $x$,  
a growth bound deduced from the contour
integral representation given in 
formula \eqn{905d} of  Appendix A implies that 
the function  $Z(\mu)= W_{\kappa, \mu}(e^{u_0})$ is
an entire function of order at most $1$.
The fact that $Z(\mu)$  is entire of order $1$ and 
maximal type will follow from the zero counting formula \eqn{513},
asserting that for some $C>0$ it has at least $C T\log T$ zeros in the
disk $|z|<T$  as $T \to \infty$.
For any entire function of order less than  $1$ or 
of order $1$ and finite type necessarily has $O(T)$
zeros in each disk $|z|\le T$ as $T \to \infty$, (cf. Titchmarsh \cite[8.75]{Ti39}).
We now establish (1)-(3). \\

(1) Apply Theorem~\ref{th21} for the Morse potential
$V_k(u)$ with $k= - \kappa$. Then the zeros of  
 $Z(\mu):=W_{\kappa, \mu}(e^{u_0})$ correspond to the $L^2$-eigenfunction
 $\psi(u_0, -\mu^2)=0$, corresponding to
  $E=-\mu^2$ being a  Dirichlet
boundary condition eigenvalue of the Morse potential
Schr\"{o}dinger equation on the 
half-line $[u_0, \infty)$. This boundary value problem
is self-adjoint, so $E$ is real. Thus $\mu$ must be real or 
pure imaginary. 

 The Schr\"{o}dinger operator spectrum is simple,
therefore each value of $E$ occurs once. Both roots in $\mu$
of $E=-\mu^2$ produce a zero of the Whittaker function,
due to  the symmetry $W_{\kappa, \mu}(x)= W_{\kappa, -\mu}(x)$.
These zeros must therefore be simple, except for 
$\mu=0$, where a double zero must occur since 
$W_{\kappa, -\mu}(x)$
is an even function of $\mu$. \\

(2) Theorem~\ref{th21} establishes that
 the Schr\"{o}dinger spectrum for Dirichlet boundary conditions is pure
discrete and is bounded below by $E = -\mu^2> -k^2$, and
by $E>0$ when $k \ge 0$.
Thus there can only be finitely many roots $\mu$ on the real axis, all
with $0 \le \mu^2 < k^2$, so that $-|\kappa| < \mu < |\kappa|$,
and there are no real roots when $k \ge 0$, i.e $\kappa \le 0$. \\

(3) The asymptotic estimate \eqn{513} follows from Theorem~\ref{th34}.
A pair of zeros  $\mu=\pm it$ corresponds to a single Dirichlet eigenvalue
$E= -\mu^2=t^2$ of the Schr\"{o}dinger operator, and the condition
$|\mu|\le T$ corresponds to the bound $E < T^2$. Thus we obtain
from the eigenvalue asymptotic \eqn{361} that
$$
N(Z(\mu)) = N(T; 0, u_0) = \frac{2}{\pi} T \log T + 
\frac{2}{\pi}(2\log 2 -1 - u_0) T + O(1),
$$
as required.
$~~~\bsq$

\paragraph{Remark.}
 For the case $k=0$,  and $\Re(w) > 0$ the Whittaker function is
related to the $K$-Bessel function (MacDonald function) by 
\beql{520a}
K_{\mu}(w) = \sqrt{\frac{\pi}{2w}} W_{0, \mu}(2w).
\eeq
The $K$-Bessel function for positive real $w=e^u$ is known to 
have all its zeros on the imaginary axis. 
According to
Erd\'{e}lyi et al, \cite[Vol II, 7.13.2 (19)]{E2}
an asymptotic formula is available when $\mu=it$ is
on the imaginary axis with
$w= x=e^u$. 
This formula states that, for $t > x >0$, 
\beql{521}
K_{it}(x) = \sqrt{2\pi} (t^2-x^2)^{-1/4} e^{-\frac{1}{2} \pi t}
\left( \sin ( t\cosh^{-1} (\frac{t}{x} ) - \sqrt{t^2-x^2}+ \frac{\pi}{4}) + 
O ( \frac{1}{ \sqrt{t^2 - x^2}} )\right).
\eeq
(We have changed the  error term, which appears to have a misprint.)
These asymptotics agree with that of Theorem~\ref{th51}.

%
%
%

\section{Eigenvalue Interpretations for Riemann Zeta Zeros}

There is considerable circumstantial evidence  for the  existence
of a spectral interpretation of the Riemann zeta zeros, see 
Berry and Keating
\cite{BK99} and Katz and Sarnak \cite{KS99}.
Indeed there  are  several known operator-theoretic interpretations  of
zeta zeros. 
There is a  scattering theory interpretation of the Riemann 
zeta zeros arising from work of 
Pavlov and Faddeev \cite{FP72} concerning the Laplacian
acting on the modular surface. A slightly
different version of this interpretation is given by Lax and Phillips \cite{LP76}, 
see in particular  \cite[Appendix 2 to
Sect. 7]{LP76}. Recently an operator formulation of the Lax-Phillips 
approach was given
by Y. Uetake \cite{Ue08}.  
There are also interpretations of the zeta zeros
in terms of eigenvalues of operators on various function spaces,
given by Connes \cite{Co99}.
The formulation of Connes permits
the construction of such an operator on a Hilbert space that  has eigenvalues
that detects only  those zeta 
zeros that are on the critical line, up to a fixed multiplicity. This approach was extended to 
automorphic $L$-functions by Soul\'{e} \cite{So99} and Deitmar \cite{De01}. 
R. Meyer (\cite{Mey05}, \cite{Mey05b}) gave an unconditional formulation of an operator on
a more general Banach space whose eigenvalues detect all zeta zeros,
including those that are off the critical line if the Riemann hypothesis fails. 
There is also recent work of Sierra (\cite{Si05}- \cite{Si08b})
and Sierra and Townsend  \cite{ST08}
concerning quantum mechanical models for the zeta zeros,
and possible ``Hilbert-Polya" operators. Burnol \cite{Bu08} 
has constructed   a specific
``toy model" operator of the integral-differential form considered by Sierra,
based on his earlier work \cite{Bu06}.

We recently observed   (\cite{La06}, \cite{La07})
that there is a natural candidate for a ``Hilbert-Polya" 
operator, using the framework of the de Branges Hilbert spaces of
entire functions, 
provided that the Riemann hypothesis
holds. This interpretation leads
to a possible connection with Schr\"{o}dinger operators on a half line,
described below.

We first review  some aspects of the de Branges theory.
The de Branges theory of Hilbert spaces of entire functions (\cite{deB68}) assigns to an entire 
function $E(z)$ of a special type  called a {\em structure function} 
a reproducing kernel Hilbert space of entire functions $\sH(E(z))$, together with  a
 (generally unbounded) linear 
operator $M_z: f(z) \mapsto zf(z)$ defined on the
 domain $\sD_z: = \{ f(z) \in \sH(E(z)): ~ zf(z) \in \sH(E(z))\}$.
This operator is closed and symmetric, with deficiency indices $(1,1)$, so that it 
has a  one-parameter family of self-adjoint extensions,
parametrized by the unit circle $S^1= U(1)$.   A {\em structure function}  $E(z)$ is any 
 entire function satisfying 
\beql{601}
|E(z)| > |E(\bar{z})|, ~~~\mbox{when}~~~ Im(z) >0.
\eeq
Functions $E(z)$ with this property are sometimes called {\em Hermite-Biehler functions}
(see Levin \cite[Chap. VII]{Le80}). Any entire function
can be uniquely represented as $E(z)= A(z) -i B(z)$, in which 
$A(z), B(z)$ are entire functions that are real on the real axis. If we
let $E^{\sharp}(z) := \overline{E(\bar{z})}$, then
$$
A(z) = \frac{1}{2} \left(E(z) + E^{\sharp}(z)\right),~~~~~B(z) = \frac{1}{2i}\left(E(z) - E^{\sharp}(z)\right).
$$
The Hermite-Biehler property \eqn{601} has the consequence that the
associated $A(z)$ and $B(z)$ have only real zeros, and these zeros interlace
(allowing multiple zeros, counting multiplicity)(cf. \cite[Lemma 2.2]{La05b}).  
The Hilbert space $\sH(E(z))$
is always nonempty, and it consists of all entire functions satisfying certain
growth restrictions (compared to $E(z), E^{\sharp}(z))$ in the upper half plane.
One particular self-adjoint extension of $(M_z, \sD_z)$ has pure discrete simple spectrum
located at the zeros of $A(z)$, and another has pure discrete simple
spectrum located at the zeros of $B(z)$. 

The de Branges theory also gives a Fourier-like transform, depending on
the structure function,  that gives an isometry to a new Hilbert space on which
the multiplication operator becomes a $2 \times 2$ matrix system of linear differential
operators, called a canonical system; we  call this the   de Branges transform.
Here  a {\em canonical system} is a 
family of differential equations, depending on $z \in \CC$ as
a parameter, on an interval $[a,b]$ in which we allow $0\le a < b \le \infty, $ given by
\beql{602}
\frac{d}{dt}\left[ {{A(t, z)}\atop{B(t,z)}} \right] =
z \bJ
 \bM(t) \left[ {{A(t, z)}\atop{B(t,z)}} \right],   
\eeq
in which 
$$
\bJ=\left[ \begin{array}{cc} 
0& -1\\ 1&0 \end{array} \right],~~~~~~
\bM(t) = \left[ \begin{array}{cc} 
\alpha(t) & \beta(t) \\ \beta(t) & \gamma(t) \end{array} \right],
$$
with $\bM(t)$ being a positive semidefinite real matrix-valued function on the interval $(a,b)$,
with coefficients being measurable real valued functions of $t$. 
In fact   de Branges's  formulation of his  theory  \cite{deB68} 
formulates results using an integral equation
which is an integrated form of the canonical system, in order to deal with
smoothness issues; he does not use canonical systems. 
However the  use of a 
canonical system makes comparison with Schr\"{o}dinger operator formulation
easier, cf. Remling \cite{Re02}. 
 In general the de Branges transform  to the canonical system is not known explicitly,
but it has been determined in a  number of "exactly solvable" examples, 
see de Branges \cite{deB64},  \cite[Chap. 3]{deB68}. 

A special subclass of canonical systems
can be  nonlinearly transformed to a pair of Schr\"{o}dinger operators on 
an interval or a half-line, a transformation
that squares the eigenvalues.  Such  canonical systems  are
a  subclass of those whose
canonical matrix $\bM(t)$ is diagonal and invertible almost everywhere.
de Branges \cite[Section 4]{deB68} proves a
 result stating that a structure function has an associated canonical
  system that is diagonal whenever $A(z)$ is an even function and
  $B(z)$ is an odd function. 
For such diagonal canonical systems,
one may 
 monotonically rescale  the time variable so that  the coefficient matrix becomes
\beql{515}
\bM(t) = \left[ \begin{array}{cc} 
\alpha(t) & 0 \\ 0 & \gamma(t) \end{array} \right],
\eeq
with determinant $1$ almost everywhere, on a rescaled interval $[\tilde{a},\tilde{b}]$ with 
$-\infty< \tilde{a}< \tilde{b} \le \infty$.
If these coefficients are smooth enough, then the transform 
gives a (supersymmetric) pair of transformed Schr\"{o}dinger operators
$- \frac{d^2}{dt^2} + V^{\pm}(t)$
having  potentials
\beql{516b}
V^{\pm}(t) := W(t)^2 \pm W'(t)=
  \frac{1}{4} \left( \frac{\alpha'(t)}{\alpha(t)}\right)^2
 \pm \frac{1}{2}\left( \frac{ \alpha''(t) \alpha(t) - \alpha'(t)^2}{\alpha(t)^2}\right)
\eeq
with associated superpotential 
$$
W(t) := \frac{1}{2}\frac{\alpha'(t)}{\alpha(t)},
$$
where $\alpha'(t) = \frac{d}{dt} \alpha(t)$.
These two Schr\"{o}dinger operators have  related boundary conditions at the left endpoint,
such that  they are self-adjoint with  identical discrete spectrum, with the exception of
the eigenvalue $0$. One of them has Dirichlet boundary conditions, the other
has boundary conditions that depend on the eigenvalue.
There are a number of known examples where, starting from
the structure function, all these steps can be carried out explicitly 
 to obtain such a pair  of Schr\"{o}dinger operators.  This is the case
 for the Morse potential for certain parameter values, 
 using the fact that it is a shape-invariant potential treatable
 by supersymmetry, cf. Cooper et al. \cite[Sec. 4.1]{CKU95}, see \cite{La08}.

In \cite{La05b},  \cite{La06}, \cite{La07} we
 considered de Branges spaces associated to the Riemann zeta function
and more general $L$-functions.    
In particular, we associate to   the Riemann $\xi$-function  a  family 
of entire  functions depending on the real parameter $h$, given for $h \ne 0$ by
$$
E_h(z) := \frac{1}{2}\left( \xi(\frac{1}{2}+h-iz) + \xi(\frac{1}{2} -h-iz)\right) - 
\frac{i}{2h}\left( \xi(\frac{1}{2}+h-iz) - \xi(\frac{1}{2} -h-iz)\right),
$$
and for $h=0$ by
\beql{604}
E_0(z):= \xi( \frac{1}{2} -iz) + \xi^{'}( \frac{1}{2} -iz),
\eeq
with $\xi^{'}(s) = \frac{d}{ds} \xi(s)$.
We show in (\cite[Theorem  2.1]{La05b}) that 
for $|h| \ge \frac{1}{2}$ each member of this family  is  unconditionally a 
de Branges structure function,  and, conditionally on the
Riemann hypothesis, that they are structure functions  for all nonzero real $h$.
Furthermore  we show in (\cite[Theorem 1]{La06}) that
the Riemann hypothesis  holds if and only if
the entire function $E_0(z)$ in \eqn{604}
is a  de Branges structure function. Assuming RH,
a theorem of de Branges (\cite[Section 47] {deB68}) applies 
to assert that the
canonical system associated to $E_0(z)$  by the de Branges transform has
matrix function $\bM(t)$ of diagonal form. 

To summarize, if the Riemann hypothesis holds, the de Branges
theory predicts the existence of a ``Hilbert-Polya" operator given
as a diagonal canonical system, specifically associated to 
the structure function \eqn{604}. Now we can ask if one 
can  nonlinearly transform this diagonal canonical system to
a pair of Schr\"{o}dinger operators. 
As obstacles, we  do not know whether  the resulting diagonal matrix function $\bM(t)$
will be almost everywhere invertible, and even if it is, we do not know whether the resulting 
matrix coefficients \eqn{515} will be smooth enough to permit the
nonlinear transform \eqn{516}. Here we speculate that the smoothness
conditions will be the main obstacle, and that one may wish to consider
this transform in the sense of generalized functions.   
(This is a suitable topic for further research in its own right.) In any case, if this transform
were possible,  there would then  exist a ``natural" 
Schr\"{o}dinger operator on a half-line,
having a  (possibly distributional) ``Hilbert-Polya" potential that, for
Dirichlet boundary conditions, has as simple eigenvalues  $\gamma^2$ 
corresponding to  $\rho=\frac{1}{2} \pm i\gamma$
being  a pair of zeta zeros.  One can then attempt to reverse-engineer this
potential, or in any case, to reverse-engineer the matrix coefficient functions of
the associated canonical system.

For  the Morse potential treated here
 the ``toy model" aspect can be extended further, to the level
  of de Branges spaces and canonical systems.
 That is, there is a family of de Branges spaces 
  associated to the Morse
potentials,  described in \cite{La08}, where the nonlinear transformation  to
a Schr\"{o}dinger operator on a half-line can be justified.
As an example, one can show that for real $u$ 
  the function
\beql{603}
E_{0}(u,z) := e^{\frac{1}{2}(e^u-u)}W_{\frac{1}{2}, i z}(e^u) 
-i ze^{-\frac{1}{2}(e^u+u)}W_{-\frac{1}{2}, i z}(e^u)
\eeq
is a de Branges structure function, and the transforms above convert it 
 to a Schr\"{o}dinger operator with
Morse potential with parameter $k= -\frac{1}{2}$ on a the corresponding half-line $[u, \infty).$ 
In particular, the self-adjoint extension of the de Branges multiplication
operator associated to
$$
A(z) = e^{-\frac{1}{2}(e^u-u)} W_{\frac{1}{2}, i z}(e^u) 
$$
is transformed  to this Schr\"{o}dinger operator with Dirichlet
boundary conditions at the left endpoint.  
This example shows that 
there is no apparent obstacle
at the level of asymptotic distribution of eigenvalues  to finding a diagonal canonical system
transformable to an associated Schr\"{o}dinger potential that encodes the
zeta zeros. Such an 
(integrated)  potential
must increase exponentially as $t \to \infty$ to have the correct asymptotics
of zero densities.

       We next address  the feature  that  these ``toy models"    reproduce only the main term
asymptotics of the distribution of zeta zeros, but not the fine structure of these zeros.
Namely, the  Riemann zeta zeros are believed to have local density statistics with fluctuations 
described by  the GUE distribution of random matrix theory, see Katz and Sarnak \cite{KS99}.
These GUE  fluctuations  encode properties of the distribution of primes,
and they partially manifest themselves in the fact that the error term  $O(\log T)$
for zeta zero asymptotics  in \eqn{112}
is nontrivial.  In contrast the 
 zeros of the Whittaker functions have normalized  local density statistics
which provably have no fluctuations, i.e. the normalized zeros will be 
completely regularly spaced, a fact which
corresponds to the $O(1)$ error term in \eqn{111a}. This fact is not so far removed from
the zeta function itself, because in \cite{La06} a similar phenomenon was observed for
zeros of differenced Dirichlet $L$-functions, whenever one moves off the critical line to a line
$Re(s) = \frac{1}{2}+ \delta$, for some  $\delta>0$. The zeros of these differenced $L$-functions
 correspond to vanishing of the real part of the $L$-function $L(\frac{1}{2} + \delta +it, \chi)$. 
The asymptotics of these differenced $L$-function zeros is known to have the same 
main term as the zeta zero asymptotics,  but also to have a $O(1)$ remainder term, just as
in these ``toy models". Furthermore \cite{La06} showed that 
these differenced $L$-function zeros do  have an eigenvalue
interpretation for a multiplication operator in a suitable de Branges space, assuming RH. 
This fact shows that the GUE phenomenon for the zeta zeros is (in some sense) confined to 
the critical line $Re(s)= \frac{1}{2}.$ We take this to mean that, assuming RH  holds, the 
de Branges space with structure function \eqn{604} will have some unusual features
of its associated canonical system. Perhaps these will be
 badly behaved coefficients, so  that if there is an associated 
``Hilbert-Polya" Schr\"{o}dinger
operator, its potential will be very badly behaved, and cannot be treated as a function. 

In any case a desire to explain  GUE  suggests further consideration of 
 the  question: what conditions must a Schr\"{o}dinger
 potential $V(x)$ satisfy in order  to have eigenvalues that
reproduce the Weyl asymptotics of the zeta zeros and  also obey local normalized
GUE statistics?

%
%
%

\section{Concluding remarks}

From the perspective of de Branges spaces, the natural operator to
consider is a canonical system. Indeed, if the Riemann hypothesis
 holds, a canonical system  ``Hilbert-Polya" operator must exist
with structure function \eqn{604}  whose (real) eigenvalues are the
imaginary parts of the zeta zeros. This paper is motivated  by the
observation that under  some conditions a canonical system can be
(nonlinearly)  transformed to
the familiar form of a one-dimensional Schr\"{o}dinger operator on
a half-line, also known as a singular Sturm-Liouville problem. 
This raises the possibility that the ``Hilbert-Polya" operator proposed for zeta zeros
might be realizable  by such a Schr\"{o}dinger  operator, which falls  squarely in the 
framework of quantum mechanics.

The existence of this nonlinear transformation
 reflects the fact that  the de Branges theory generalizes 
part of Weyl's theory of singular Sturm-Liouvlille problems, cf.
Everitt and Kalf \cite{EK08}. The Morse potential example considered
here is  a singular Sturm-Liouville problem which falls within
Weyl's theory. That theory assigns to each such problem an analytic
function, the {\em principal Weyl-Titchmarsh $m$-function}, which encodes information about
the spectral data of the problem, encoded in the initial conditions.
This function is discussed in Appendix B. This function  is analytic in the upper half-plane
and has positive imaginary part there, i.e. it is a so-called Herglotz function.

One may also associate  directly to a de Branges structure function $E(z)$  the function 
\beql{605}
m_{E}(z) := - \frac{B(z)}{A(z)},
\eeq
which we will here call  the {\em de Branges $m$-function}.
It is a Herglotz function associated to 
the de Branges multiplication operator
$(M_z, \sD_z)$  which plays a role analogous to that of the principal Weyl-Titchmarsh
$m$-function for Schr\"{o}dinger operators. That is,  it encodes certain  information on the spectrum 
of self-adjoint extensions of this symmetric operator. More precisely, under the de Branges transform
it equals a ``classical"  $m$-function attached to a Dirac operator,  in the sense of
Levitan and Sargsan \cite[Chap. 3]{LS75}. This Dirac operator  is directly constructed
from the canonical system   \eqn{602}, by a recipe we omit here. In 
the case of the de Branges structure function \eqn{603} above, using a
Whittaker function identity \cite[(13.4.30)]{AS66}, we obtain 
\beql{606}
-\frac{B(z)}{A(z)} = -e^{-e^u} \frac{zW_{-\frac{1}{2}, iz}(e^u)}{W_{\frac{1}{2}, iz}(e^u)}
=  - \frac{1}{z} \left[ e^{-e^u}( -\frac{W_{\frac{3}{2}, iz}(e^u)}{W_{\frac{1}{2}, iz}(e^u)} + e^u-1 )\right]
\eeq
For comparison, in Appendix B we  determine the principal 
 Weyl-Titchmarsh $m$-function for the Morse potential 
 on the right half-line (Theorem~\ref{th31}). The formula \eqn{606} should
 be compared with Theorem~\ref{th31} 
with parameter $k=-\frac{1}{2}$. 
The computations in Appendix B  
may be useful as a guide for determining  relations
between Weyl-Titchmarsh $m$-functions 
and suitable de Branges $m$-functions.


%
%
%

\section{Appendix A: Whittaker Functions in the Complex Domain}

We recall here facts about  Whittaker functions $W_{\kappa, \mu}(z)$ as
functions of three complex variables $(\kappa, \mu,z)$.
They are a type of 
confluent hypergeometric
function, and are single-valued in the variables $\kappa$
and $\mu$, and
multi-valued in the $z$ variable, with
a singular point at $z=0$.

The standard confluent hypergeometric function
${}_1 F_1 (\alpha, \beta; z)$ is given by the power series expansion
\beql{901}
{}_1 F_1 (\alpha, \beta; z) := \frac{\Gamma(\beta)}{\Gamma(\alpha)}
\sum_{j=0}^{\infty} \frac{\Gamma(\alpha +j)}{\Gamma(\beta+j)} \frac{z^j}{j!}.
\eeq
The  normalized function
\beql{902}
{}_1 \tilde{F}_1 (\alpha, \beta; z) := 
\frac{1}{\Gamma(\beta)} {}_1 F_1 (\alpha, \beta; z) 
\eeq
is  an entire function of three complex variables
$(\alpha, \beta, z) \in \CC^3$,
as can be read off from the uniform convergence
properties of the
expansion \eqn{901} on compact subsets of
$\CC^3$, observing  that $\frac{1}{\Gamma(\beta+j)}$ is entire and 
$\frac{\Gamma(\alpha +j)}{\Gamma(\alpha)}$
is  a  polynomial in $\alpha$,  cf. Buchholz \cite[Sec. 1.3]{Bu69}.
(The unnormalized function ${}_1 F_1 (\alpha, \beta; z)$ has simple
poles at $\beta=0, -1, -2, ...$ for most  $(\alpha, z)$.)
It satisfies the 
confluent hypergeometric differential equation
\beql{903}
z\frac{d^2F}{dz^2}+ (\beta -z) \frac{dF}{dz} - \alpha F =0.
\eeq 
which has a regular singular point at $z=0$ and an irregular
singular point at $z=\infty$, of irregularity index one. (Slater \cite{Sl60}). 

The Whittaker functions  
$M_{\kappa, \mu}(z)$ and  $W_{\kappa, \mu}(z)$ were given
by Whittaker \cite{Wh04}  in 1904 as particular solutions to
the differential equation
\beql{904}
\frac{d^2F}{dz^2} + \left( - \frac{1}{4} + \frac{\kappa}{z} +
 \frac{\frac{1}{4} - \mu^2}{z^2} \right) F = 0.
\eeq
This equation also  has  a
regular singular point at $z=0$ and an irregular singular point
at $z= \infty$ on the Riemann sphere, of irregularity index one.
Whittaker sets 
\beql{905}
M_{\kappa, \mu}(z) := e^{-\frac{z}{2}} z^{\frac{1}{2}+\mu} 
{}_1 F_{1} (\frac{1}{2} -\kappa +\mu, 2\mu+1; z),
\eeq
with $z^{-\mu} = e^{\mu \log z}$, and 
the $z$-plane is cut along the negative real axis.
This defines the {\it principal branch} of the
Whittaker $M$-function. 
Under analytic continuation in the
$z$-variable, this  function is multivalued,
with multivaluedness arising entirely from the function $ z^{-\mu}$.
It also has poles at $\mu=0, -\frac{1}{2}, -1, ...$ and to 
eliminate these (excluding $\mu=0$)
Buchholz \cite[p. 12]{Bu69} introduces the
normalized  function $\sM_{\kappa, \mu}(z)$, defined by
\beql{905b}
\sM_{\kappa, \mu}(z) := \frac{1}{\Gamma(1+2\mu)}M_{\kappa, \mu}(z).
\eeq

The  Whittaker
function $W_{\kappa, \mu}(z)$ 
 is specified by the asymptotic property of having rapid decrease as
$z=x \to \infty$ along the positive real axis. Whittaker \cite{Wh04}
defined it using two integral representations given in
Whittaker and Watson \cite[Sec. 16.12]{WW63}. In terms of
the confluent hypergeometric
function above we have (Truesdell \cite[p. 170]{Tr48})
\begin{eqnarray}~\label{906}
W_{\kappa, \mu}(z) &= & 
\frac{\Gamma(-2\mu)z^{\frac{1}{2} +\mu} e^{- \frac{z}{2}}}
{\Gamma( \frac{1}{2} -\kappa-\mu)} 
{}_1 F_1 ( \frac{1}{2} -\kappa+\mu, 2\mu+1; z) \nonumber \\
&& 
~~~~ + \frac{\Gamma(2\mu)z^{\frac{1}{2} -\mu} e^{- \frac{z}{2}}}
{\Gamma( \frac{1}{2} -\kappa+\mu)} 
{}_1 F_1 ( \frac{1}{2} -\kappa-\mu, -2\mu+1; z).
\end{eqnarray}
This formula can alternatively be written 
\begin{equation}~\label{906aa}
W_{\kappa, \mu}(z) =  
\frac{\Gamma(-2\mu)}{\Gamma(\frac{1}{2} -\kappa -\mu)}M_{\kappa, \mu}(z)+ 
\frac{\Gamma(2\mu)}{\Gamma(\frac{1}{2} -\kappa+\mu)}M_{\kappa, -\mu}(z),
\end{equation}
which incidentally shows that $W_{\kappa, \mu}(z) = W_{\kappa, -\mu}(z)$.
The formula \eqn{906} directly defines 
$W_{\kappa, \mu}(z)$ as an analytic  function of
three complex variables,  in the region $(\kappa, \mu) \in \CC^2$
and all $z$ in the complex plane cut along the negative
real axis. ( The analyticity
is read off using the function \eqn{902}.)
We term $W_{\kappa,\mu}(z)$ on this region the
{\em principal branch} of the Whittaker function.
Under analytic continuation of this function
in the $z$-variable, this
function is multivalued, and becomes
single-valued on a logarithmic covering surface around
the point $z=0$, remaining an entire function
of $\kappa$ and $\mu$ on every branch above
a fixed $z$. (For certain  values of $\kappa$ and $\mu$ it may
be single-valued on a finite
 cover of $\CC^{\ast} = \CC\backslash \{0\}$.)

The function $W_{\kappa, \mu}(z)$ also 
has the contour integral representation
(\cite[Sec. 16.12]{WW63})
\beql{905d}
W_{\kappa, \mu}(z) = \Gamma(\frac{1}{2}+\kappa - \mu) e^{-\frac{z}{2}}
z^{\kappa} \left( -\frac{1}{2\pi i} 
\int_C (-t)^{-\kappa- \frac{1}{2} +\mu}
(1 + \frac{t}{z})^{\kappa - \frac{1}{2} +\mu}e^{-t}dt\right),
\eeq
in which the contour $C$ goes from $+\infty$ to $0^{+}$ and
back by a ``keyhole'' contour around the positive real axis,
keeping $t=-z$ outside the contour, with $arg(z)$ taking the
principal value, $|arg(-t)|\le \pi$, and $arg(1+ \frac{t}{z})$
taking that value which goes to zero at $t \to 0$ by
a path outside the contour. If we restrict $z$ to the positive
 real axis, one obtains by estimating the
contour  integral  \eqn{905d} that
$$
W_{\kappa, \mu}(z)= O \left(e^{C |\mu |\log (|\mu|+2)}\right),
$$
viewed as a function of $\mu$, for fixed $\kappa$ and $z$, 
a constant $C$ depending on $\kappa$ and $z$. This implies 
that $W_{\kappa, \mu}(z)$ is an entire function of order at most $1$
in the $\mu$-variable.

The distinguishing property of
the principal branch of $W_{\kappa, \mu}(z)$,
  is its rapid  decrease
 along the real axis $z=x$ as $x \to \infty$.  All other linearly independent
solutions to Whittaker's equation \eqn{904}  increase rapidly 
 (in absolute value) along the positive real axis as $x \to \infty$.
The rapid decrease
property of the principal
branch of the Whittaker function holds 
uniformly   as $|z| \to \infty$   on the
angular  sector 
$|\arg(z)| < \frac{\pi}{2} - \epsilon$ for
any fixed $\epsilon >0$.
Furthermore, for fixed $(\kappa,\mu)$, it has on the larger angular 
sector 
$-\pi + \epsilon < \arg(z)< \pi - \epsilon$,
 an asymptotic expansion (\cite[Sec. 16.3]{WW63})
which gives
\beql{908}
W_{\kappa ,\mu}(z) = e^{- \frac{z}{2}} z^{\kappa}\left( 1 + O \left( 
\frac{\mu^2- (\kappa- \frac{1}{2})^2}{z}\right) \right).
\eeq
The Whittaker function solutions having  this
rapid decrease property glue together smoothly  in the $(\kappa,\mu)$
variables  (uniformly with $|\kappa|^2 + |\mu|^2 < R$ for any fixed $R$)
to form entire functions of $(\kappa, \mu) \in \CC^2$, 
yielding  the principal branch defined
above.
For certain values of $\kappa$ and $\mu$ we have
a finite covering of $\CC^{\ast}$;
in these cases there are a finite number of other branches. \\

%
%
%

\section{Appendix B: Weyl-Titchmarsh m-Functions for Morse Potentials }

For Schr\"{o}dinger operators on the half-line $[u_0, \infty)$
with a real-valued potential, the
Weyl-Titchmarsh $m$-function $m(u_0,E)$ is an analytic function
defined on $E \in \CC \backslash \RR$, which 
encodes information about the
spectrum of the Schr\"{o}dinger operator for various
constant boundary conditions. 
It maps the upper half-plane  $\Im(E) >0$ into itself,
and its values on the lower-half plane $Im(E) <0$
satisfy $m(u_0, \overline{E}) = \overline{m(u_0, E)}$.
Spectral information is extractible from its asymptotics 
in the upper half plane as the
real axis is approached, see Gilbert and Pearson \cite{GP87}.
In general it is not defined on the real axis, but 
in the case where the spectrum is pure discrete then
it has a meromorphic extension to the whole plane $\CC$,
which is real-valued on the real axis.

To define it, 
let $\theta_{\alpha}(u; E), \varphi_{\alpha}(u, E)$ denote the
solutions to the initial value problems
$$
\left( - \frac{d^2}{du^2} + V(u)\right) \psi(u) = E \psi(u),
$$
with
$$
\left[ \begin{array}{c}
\theta_{\alpha}(u_0, E) \\
\theta_{\alpha}^{'}(u_0, E) \end{array} \right] =
\left[  \begin{array}{c}
\cos\alpha \\
\sin\alpha  \end{array} \right]
 ~~\mbox{and}~~
\left[ \begin{array}{c}
\varphi_{\alpha}(u_0, E) \\
\varphi_{\alpha}^{'}(u_0, E) \end{array} \right] =
\left[  \begin{array}{c}
-\sin \alpha \\
\cos \alpha  \end{array} \right].
$$
We define  the $m$-function $m_{\alpha}(E):=m_{\alpha}(u_0,E)$ 
by the condition  that
\beql{330}
\psi(u_0, E) := \theta_{\alpha}(u_0, E)+ 
m_{\alpha}(E) \varphi_{\alpha}(u_0, E)
\eeq
belongs to $L^2([u_0, \infty);du); $ there is a unique solution when
$Im(E) \ne 0$.  The case $\alpha=0$ 
gives the {\em principal Weyl-Titchmarsh $m$-function},
where  $\varphi_{0}$ satisfies Dirichlet boundary
conditions at the left endpoint,
and $\theta_{0}$ satisfies Neumann boundary conditions;
we write $m(u_0,E) := m_{0}(u_0,E)$. For the principal
$m$-function we have the alternate formula
\beql{331}
m(u_0, E) = \frac{\psi'(u_0, E)}{\psi(u_0, E)} 
= \frac{d}{du} \log \psi(u, E) \Bigr|_{u=u_0},
\eeq
which follows since $\psi(u_0, E) =1$ and $\psi'(u_0, E) = m(u_0, E)$.
The different $m$-functions for different $\alpha$ are related
by linear fractional transformations, namely
\beql{332}
m_{\alpha}(u_0, E) = \frac{ (\cos \alpha) m_0(u_0, E) - \sin \alpha}
{(\sin \alpha) m_0(u_0, E) + \cos \alpha},
\eeq
the principal $m$-function 
determines them all. With our definition the $m$-function
maps the upper half-plane $\CC^{+} =\{ E: \Im(E) >0\}$
into itself, so that it is a Herglotz function.  (Our definition
is the negative of the $m$-function for $\alpha$ given
in Titchmarsh \cite[Sec. 2.1]{Ti62}, which maps the upper half-plane
to the lower half-plane $\CC^{-}$, cf. \cite[Sec. 3.1]{Ti62}.)
As noted above, in  the general case an $m$-function need not analytically continue 
across the real axis anywhere. 
However in  the case at hand, that of   pure discrete spectrum, 
 the Weyl-Titchmarsh $m$-function extends to  a meromorphic
function on $\CC$, and all of  its singularities are simple poles, located
on the real axis. This case is 
treated in Chapter II of Titchmarsh \cite{Ti62}. 
The principal $m$-function $m(u_0,E)$ then has a pole at any
point where $\psi(u_0, E)$ would be a multiple of $\varphi(u_0, E)$.

%
%
%

\begin{theorem}~\label{th31}
The Schr\"{o}dinger equation on the half-line
$[u_0, \infty)$ with Morse potential $V_k(u) = \frac{1}{4} e^{2u}+ ke^{u}$
has principal Weyl-Titchmarsh  $m$-function $m(u_0, E)$ with
eigenvalue $E= z^2$ for  endpoint $u_0$ given by 
\beql{333}
m(u_0, z^2) = 
-\frac{ W_{1-k, iz}(e^{u_{0}}) }{W_{-k, iz}(e^{u_0})}
+ \left(\frac{1}{2} e^{u_{0}}+ k  -\frac{1}{2} \right)
\eeq
\end{theorem}

\paragraph{Proof.}
Since the Morse potential on the half-line $[u_0, \infty)$ has pure discrete spectrum, the
$m$-function will be a meromorphic function of the variable
$E=z^2$. 
Theorem ~\ref{th21}(1) shows that for $\Im(E) > 0$
the subordinate solution is given by 
$$
\psi(u, z^2) := e^{-\frac{u}{2}} W_{-k, iz}(e^{u}) = 
e^{-\frac{u}{2}} W_{-k, -iz}(e^{u}). 
$$
We conclude that 
\begin{eqnarray}~\label{334}
m(u, z^2) & =&  \frac{\psi^{'}(u, z^2)}{\psi(u, z^2)} \nonumber\\
&=& 
- \frac{1}{2} + e^u\frac{ W_{-k, iz}^{'}(e^{u})}{ W_{-k, iz}(e^u)},
\end{eqnarray}
in which $W_{\kappa, iz}^{'}(x) = \frac{d}{dx}   W_{\kappa, iz}(x).$
We now use the identity (\cite[13.4.33]{AS66}) 
\beql{335}
xW_{\kappa, \mu}^{'}(x) = (\frac{1}{2}x - \kappa)W_{\kappa, \mu}(x) 
- W_{\kappa+1, \mu}(x).
\eeq
Taking $x=e^u$ and substituting  this in \eqn{334} results in
\beql{336} 
m(u, z^2)= \left(\frac{1}{2}e^u+k  - \frac{1}{2}\right)
- \frac{ W_{1-iz, \mu}(e^u)} {W_{-k, iz}(e^u)},
\eeq
the desired formula.$~~~\bsq$

\paragraph{Remarks.}
(1) The principal $m$-function contains enough information to
uniquely reconstruct the potential $V(u)$ on $[u_0, \infty)$, under
very general conditions; this is a particular kind of inverse spectral problem.
In the case where the Schr\"{o}dinger
operator has  pure discrete spectrum, this follows from 
the fact that the Dirichlet problem eigenvalues
are the poles of $m(u_0, E)$, and the Neumann problem
eigenvalues are the zeros of $m(u_0, E)$; necessarily all
of these fall on the real axis. The Dirichlet and
Neumann spectrum together are known to be enough
to determine the spectrum. \\

(2) An important feature of the principal $m$-function
is that it satisfies a Riccati equation in the $u$-variable,
\beql{337}
\frac{d}{du} m(u, E) + m(u, E)^2 = V(u) - E.
\eeq
This can be exploited in inverse spectral reconstruction
of the potential. \\

(3) It is well-known that the principal $m$-function
attached to a Schr\"{o}dinger
operator on a half-line having a  potential $V(u)$ that 
is locally $L^1$ (and which is in the limit point
case at $\infty$) has strong restrictions  on its
asymptotics in the upper half-plane as $\Im(E) \to \infty$.
Let $E= z^2$ with $z$
confined to a sector $0< \epsilon \le\arg{z} \le \frac{\pi}{2}- \epsilon.$ 
Then for $a >0$ and $0 \le u_0 < a$,
for sufficiently large $|z|$ with 
$\epsilon < \arg(z) < \frac{\pi}{2} - \epsilon$, there holds
\beql{338}
m(u_0, z^2) = iz  - \int_{0}^{a-u_0} V(t+u_0) e^{2izt}dt +
O\left(\frac{1}{|z|}\right),
\eeq
where the constant in the $O$-symbol depends on $a$ and $\epsilon$.
The basic result asserting $m(u_0, E)$ equals $i \sqrt{E} $
plus an error term  is due to Atkinson \cite[Theorem 2 ff.]{At81},
who deduced it using \eqn{337}. Atkinson 
 used a different definition of  $m$-function, 
which in terms of the definition used here is $\frac{-1}{m(u, z)}$.
(That is, Atkinson used the Neumann boundary condition in defining his $m$-function,
while the principal $m$-function is based on the Dirichlet boundary condition at the left
endpoint.)
A result in the form  \eqn{338}  appears in 
Rybkin \cite[Corollary 5.2]{Ry02}. 
Combining  Theorem~\ref{th31} with the asymptotics \eqn{338} 
provides information on the relative sizes  of two Whittaker
functions $W_{1-k, \mu}(e^u)$ and  $W_{-k, \mu}(e^u)$
for $\mu$ in the second quadrant region of the $\mu$-plane. \\

%
%
%

\end{document}